\documentclass[12pt]{iopart}

  \usepackage{paralist}
  \usepackage{graphics} 
  \usepackage{epsfig} 
\usepackage{graphicx}  \usepackage{epstopdf}

\usepackage[colorlinks=true]{hyperref}

\newtheorem{thm}{Theorem}[section]

\newtheorem{prop}[thm]{Proposition}

\newtheorem{theorem}{Theorem}[section]

\newtheorem{lemma}[theorem]{Lemma}

\newtheorem{conjecture}[theorem]{Conjecture}
\newtheorem{definition}[theorem]{Definition}
\newtheorem{remark}[theorem]{Remark}

\newcommand{\eps}{\varepsilon}

\def\debproof{ {\bf Proof.} }
\def\finproof{\hfill {\small $\Box$} \\}  

\def\finproofo{\hfill {\small $\Box$} \\}  
 
\usepackage{graphicx,color}
\usepackage{epstopdf}
\usepackage{pifont}

\newcommand{\eqref}[1]{(\ref {#1})}

\DeclareMathAlphabet{\itbf}{OML}{cmm}{b}{it}

\def\bq{{{\itbf q}}}

\def\br{{{\itbf r}}}

\def\by{{{\itbf y}}}
\def\bx{{{\itbf x}}}

\def\bY{{{\itbf Y}}}
\def\bX{{{\itbf X}}}

\def\bk{{{\itbf k}}}

\def\bW{{{\itbf W}}}
\def\bR{{{\itbf R}}}

\renewcommand{\d}{{\rm d}}

\newcommand{\RR}{\mathbb{R}}

\newcommand{\EE}{\mathbb{E}}
\newcommand{\cA}{\mathcal {A}}

\newcommand{\ea}{\end{eqnarray}}  
\newcommand{\ba}{\begin{eqnarray}}  
\newcommand{\ean}{\end{eqnarray*}}  
\newcommand{\ban}{\begin{eqnarray*}}  
\newcommand{\iint}{\int\int}

\usepackage{iopams}  
\begin{document}

\title[Imaging  by  speckle intensity correlations]
      {Imaging through a scattering medium by  speckle intensity correlations }

\author{Josselin Garnier$\hbox{}^{(1)}$  and  Knut S\o lna$\hbox{}^{(2)}$}

\address{$\hbox{}^{(1)}$Centre de Math\'ematiques Appliqu\'ees, Ecole Polytechnique, 91128 Palaiseau Cedex, France}
\ead{josselin.garnier@polytechnique.edu}
 \vspace{10pt}
\address{$\hbox{}^{(2)}$Department of Mathematics, 
University of California Irvine, Irvine CA 92617}
\ead{ksolna@math.uci.edu}
 \vspace{10pt}
  
 \begin{abstract}
 In this paper we analyze an imaging technique based on 
intensity speckle correlations over incident field position proposed 
 in [J. A. Newmann and K. J. Webb, Phys. Rev. Lett. 113, 263903 (2014)].
Its purpose is to reconstruct a field incident on a strongly scattering random medium. 
The  thickness of the complex medium is much larger than the 
scattering mean free path so that the wave emerging from
the random section forms an incoherent  speckle pattern.
 Our analysis  clarifies 
 the conditions under which the method can give a good reconstruction
 and characterizes   
its performance.   The analysis is carried out in the white-noise paraxial regime, 
which is relevant for the applications in optics that motivated the original paper.\end{abstract}

\vspace{2pc}
\noindent{\it Keywords}: Waves in random media, speckle imaging, multiscale analysis.

\maketitle

 \section{Introduction}

Imaging and communication through a randomly scattering medium is challenging
because the coherent incident waves are transformed into incoherent wave fluctuations.
This degrades wireless communication \cite{alamouti,book1}, medical imaging \cite{huang}, and 
astronomical imaging \cite{tokovinin}.
When scattering is weak, 
different methods have been proposed, which consists in extracting the small coherent wave from the recorded field \cite{aubry09,aubry11,borcea11,borcea05,sha14}.
These methods fail when scattering becomes strong and the coherent field completely vanishes.
However recent developments have shown that it is possible to achieve wave focusing through a strongly scattering medium 
 by control of the incident wavefront
\cite{vellekoop10,vellekoop07,vellekoop08}.
These results have opened the way to new methods for wave imaging through a strongly scattering medium \cite{katz12,mosk12,popoff10}.

In \cite{webb14} an original imaging method is presented 
that makes it possible to reconstruct fields incident on a  randomly scattering medium from intensity-only measurements.
From the experimental point of view, 
the speckle intensity images are taken as a function of incident field position 
and then used to calculate the speckle intensity correlation over incident
position.
From the theoretical point of view, 
the speckle intensity correlation function is then expressed using a moment theorem as the magnitude squared of the incident 
 field autocorrelation function. The modulus of the spatial Fourier transform of the incident field can then be extracted, 
 and the incident field itself can be reconstructed using a phase retrieval algorithm.
The key argument is the moment theorem that is based  on a zero-mean circular Gaussian assumption for the transmitted field.
In \cite{webb14} the authors claim that heavy clutter is necessary and sufficient for this.
One of the main applications is a new method to view binary stars from Earth
(using the Earth's rotation and atmospheric scatter). Other biomedical applications are proposed
and extensions of the technique to imaging hidden objects 
with speckle intensity correlations over object position
have been proposed~\cite{webb16}.

In this paper we present a detailed analysis of the technique in the white-noise paraxial regime,
which is the regime relevant for the applications \cite{strohbehn,tappert}.
We clarify the conditions under which the imaging approach proposed in \cite{webb14} can be efficient.
In particular, we will see that the zero-mean circular Gaussian assumption is not strictly necessary,
however,  that strongly scattering media may not create the right conditions for the 
imaging approach to work well. 
We can distinguish two strongly scattering regimes, the scintillation regime (in which the correlation radius
of the medium fluctuations is smaller than the field radius) and 
the spot-dancing regime (in which the correlation radius
of the medium fluctuations is larger than the field radius), and these regimes
 give completely different results.
In the scintillation regime we will explain that 
the method proposed by \cite{webb14} can give a correct image, but not in the spot-dancing regime.
In particular, the spot-dancing regime may be relevant for Earth-based astronomy \cite{andrews}, 
which would let little hope that the method can be used there, but 
it could be efficient in other configurations in the scintillation regime.

The paper is organized as follows.
In Section \ref{sec:intcor} we describe the experiment and introduce the empirical speckle intensity covariance.
In Section \ref{sec:ito} we present the white-noise paraxial wave equation.
We analyze the properties of the statistical speckle intensity covariance in the 
scintillation regime in Section \ref{sec:scin}
and in the spot-dancing regime in Section \ref{sec:spot}. 
Section \ref{sec:con} summarizes the main findings.

\begin{figure}
\begin{center}
\begin{tabular}{c}
\includegraphics[width=6.2cm]{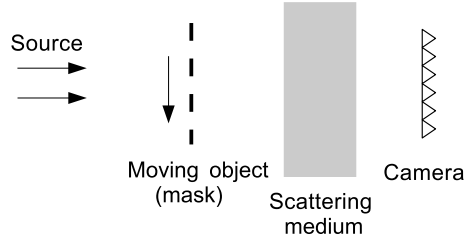}  
\end{tabular}
\end{center}
\caption{The experimental imaging set-up. The source transmits a time-harmonic plane wave.
The object to be imaged is a mask.
For each position of the mask the intensity of the transmitted field can be recorded 
by the camera. }
\label{fig:1} 
\end{figure}

\section{The intensity covariance function}
\label{sec:intcor}%
The spatial variable is denoted by $(\bx,z)\in \RR^d \times \RR$.
The source transmits a time-harmonic plane wave going into the $z$-direction with frequency $\omega$
and wavenumber $k_o=\omega/c_o$, with $c_o$ the background velocity.
The object to be imaged is a mask (a double slit in the experiment \cite{webb14}) that can be shifted transversally by a
shift vector denoted by $\br$ so that the field just after the mask is  of the form
\begin{equation}
\label{eq:inc}
U_\br (\bx) = U(\bx-\br) ,
\end{equation}
for some function $U$ (see Figure \ref{fig:1}). 
{Note that we here assume that the homogeneous scattering medium  
fills the space in between the mask and the camera, see also Remark~\ref{remark:49}. } 

The time-harmonic field in the plane of the camera is denoted by $E_\br ( \bx ) $.
It results from the propagation of the incident field $U_\br$ through the scattering medium.
The measured or empirical intensity covariance is
\begin{eqnarray}
\nonumber
C_{\br,\br'} &=& \frac{1}{|A_o|}
\int_{A_o} |E_\br ( \bx_0 )|^2 |E_{\br'}( \bx_0 )|^2 \d \bx_0
\\
&&
- \Big( \frac{1}{|A_o|}
\int_{A_o} |E_\br ( \bx_0 )|^2  \d \bx_0 \Big)
\Big( \frac{1}{|A_o|}
\int_{A_o} |E_{\br'} ( \bx_0 )|^2  \d \bx_0 \Big) ,
\label{def:intcor}
\end{eqnarray}
where $A_o$ is the spatial support of the  camera.
The conjecture found in \cite{webb14} is the following one.
\begin{conjecture}    
  \ba
\label{eq:pred}
C_{\br,\br'}  \approx    \Big| 
\int_{\RR^d}
|\hat{U}(\bk)|^2 \exp \big( i\bk \cdot ( \br'-\br) \big)  \d\bk
 \Big|^2 ,
\ea
 up to a multiplicative constant, where
\begin{equation}
\hat{U}(\bk) = \int_{\RR^d} U(\bx) \exp \big( -i \bk \cdot \bx \big) \d \bx .
\end{equation}
\end{conjecture}
When this formula holds, it is possible to reconstruct the incident field $U$
by a phase retrieval algorithm as shown in \cite{webb14}. 
Indeed (\ref{eq:pred}) gives the modulus of the inverse Fourier transform of $|\hat{U}(\bk)|^2$,
and we know the phase of $|\hat{U}(\bk)|^2$,  which  is zero, so that 
a Gerchberg-Saxon-type iterative algorithm can be applied to reconstruct $|\hat{U}(\bk)|^2$ \cite{fienup,fienup87}.
Using the estimated value of the modulus of the Fourier transform of $U(\bx)$ and applying again the same algorithm 
(assuming that the phase of $U(\bx)$ is known, for instance, equal to zero) 
it is possible to extract the incident field $U(\bx)$.
The main question we want to address is to understand under which circumstances 
and to what extent the formula (\ref{eq:pred}) holds  true.

In the expression (\ref{def:intcor}) it is assumed that the pixel size of the camera is so small that 
it is possible to consider that the camera measures the spatially resolved intensity pattern.
It is of interest to address the role of the pixel size and to assume  that the measured 
intensity is rather 
\begin{equation}
I_\br^{\rho_o}(\bx_0) =
\frac{1}{(2\pi)^{d/2} \rho_o^d} \int_{\RR^d} |E_\br ( \bx_0 +\by_0 )|^2 \exp\Big(- \frac{|\by_0|^2}{2\rho_o^2}\Big) \d \by_0 ,
\end{equation}
where $\rho_o$ is the size of the pixel of the camera.
Then the measured or empirical intensity covariance is
\ba
\label{def:Crrrho0}   \hspace*{-2.3cm}
C_{\br,\br'}^{\rho_o} = \frac{1}{|A_o|}
\int_{A_o} I_\br^{\rho_o}(\bx_0) I_{\br'}^{\rho_o}(\bx_0)\d \bx_0
- \Big( \frac{1}{|A_o|}
\int_{A_o} I_\br^{\rho_o}(\bx_0)  \d \bx_0 \Big)
\Big( \frac{1}{|A_o|}
\int_{A_o} I_{\br'}^{\rho_o}(\bx_0)  \d \bx_0 \Big)  .
\ea
Note that in order to characterize $C_{\br,\br'}^{\rho_o}$ we need to be able to evaluate
fourth-order moments for the field $E_\br ( \bx ) $. 
We describe in the next section the It\^o-Shr\"odinger model 
that makes it possible to compute such fourth-order moments and in particular how we can use this 
to characterize the   intensity covariance.  In Sections \ref{sec:scin} and \ref{sec:spot} we delineate 
two important sub-regimes of the It\^o-Shr\"odinger model corresponding respectively
to a large  or small radius of the mask,  and
show how the measured intensity covariance function can be characterized in these
cases based on our general theory for the fourth moment.

\section{The white-noise paraxial model}
\label{sec:ito}%
The model for the time-harmonic field in the plane of the camera is
\begin{equation}
\label{eq:model1}
E_\br ( \bx ) = \int_{\RR^d} \hat{g}\big((\bx,{\ell}), (\bx',0) \big) U_\br(\bx') \d \bx' ,
\end{equation}
where $U_\br$ is the incident field (\ref{eq:inc}) {in the plane $z=0$}, ${\ell}$ is the propagation distance to the camera
{localized in the plane $z=\ell$},
and $\hat{g}$ is the 
fundamental solution of the white-noise paraxial wave equation 
which we describe in the next subsections. There should be an additional factor $\exp(i k_o {\ell})$ in (\ref{eq:model1})
but it does not play any role as we only record intensities.

\subsection{The random paraxial wave equation}
We consider the time-harmonic form of the scalar wave equation {with a source of the form $2ik_o f(\bx)\delta(z)$
localized in the plane $z=0$ (which corresponds to an initial condition 
for the field of the form $f(\bx)$ in the plane $z=0$ as we will see below)}:
\begin{equation}
\label{eq:wave0}
(\partial_z^2+\Delta) E+  k_o^2 \big(1 + \mu(\bx,z)\big) E = 2ik_o \delta(z) f(\bx) ,
\end{equation}
where $\Delta$ is the transverse Laplacian (i.e., the Laplacian in $\bx$) and $f$ is a source in the plane $z=0$.
Here $\mu$ is a zero-mean, stationary, $d+1$-dimensional random process
with mixing properties in the $z$-direction {(this means that we assume that the medium is statistically homogeneous
from the plane $z=0$ to the plane $z=\ell$)}.
The function $\hat\phi $ (slowly-varying envelope of a plane wave going along the $z$-axis) defined by
\begin{equation}
E ( \bx,z)  = e^{i k_o  z  }
 \hat\phi  \big(  {\bx} ,z \big) 
\end{equation}
satisfies
\begin{equation}
\label{eq:bitos}
 \partial_{z}^2  \hat\phi+  
\left( 2 i k_o  \partial_z \hat\phi + \Delta  \hat\phi  +  k_o^2  \mu\big(  \bx  ,  z \big)  \hat\phi  \right)= 2i k_o \delta(z) f(\bx)  .
\end{equation}

\begin{definition}
\label{def:par} 
In the white-noise paraxial regime, the 
wavelength is much smaller than the initial field radius and the correlation radius of the medium,
which are themselves much smaller than the propagation distance, 
{in such a way that the product of the wavelength and the propagation distance is of the same order
as  the square radii.}
\end{definition}
In the white-noise paraxial regime, the forward-scattering approximation in direction $z$ is valid (i.e., the second derivative in $z$ in (\ref{eq:bitos}) can be neglected)
and the white-noise approximation is valid (i.e., $\mu$ can be replaced by a white noise in $z$),
so that $\hat\phi$ satisfies the 
It\^o-Schr\"odinger equation \cite{garniers1}
\begin{equation}
\label{eq:IS}
 2 i  k_o \d_z \hat{\phi} ( \bx,z) +\Delta \hat{\phi} ( \bx,z) \d z + k_o^2  \hat{\phi}( \bx,z)\circ \d B(\bx,z) =0  ,  
\end{equation}
starting from $  \hat{\phi} ( \bx,0)    = f(\bx)$, 
where $B(\bx,z)$ is a Brownian field,  that is,  a Gaussian process with mean zero and covariance function
\begin{equation}
\label{def:covgaus}
\EE\big[ B(\bx,z)B(\bx',z')\big] =  \gamma_0(\bx-\bx') \big( z \wedge z' \big),
\end{equation}
with
\begin{equation}
\label{def:gamma0}
\gamma_0(\bx)=  \int_{-\infty}^\infty  \EE[ \mu({\bf 0},0)\mu(\bx,z) ] \d z .
\end{equation}
Here the $\circ$ stands for the Stratonovich stochastic integral. 
The rigorous statement has the form of a convergence theorem
for Hilbert-space valued processes \cite{garniers1}.

\subsection{The fundamental solution}
The fundamental solution $\hat{g}$ is defined as the solution
of the It\^o-Schr\"odinger equation in $(\bx,z)$:
\begin{equation}
\label{def:greens}
 2i k_o \d_z \hat{g} + \Delta \hat{g}  \d z+ k_o^2 \hat{g} \circ \d B(\bx,z)= 0, 
\end{equation}
starting from $\hat{g}\big( (\bx,z=z'),(\bx',z' ) \big) = \delta(\bx-\bx')$.
In a homogeneous medium ($B \equiv 0$) the fundamental solution is (for $z> z'$)
\begin{equation}
\label{eq:green0}
\hat{g}_0 \big(  (\bx,z), (\bx',z') \big) =
\Big( \frac{ k_o}{2 i \pi (z-z')} \Big)^{d/2}  \exp \Big(  i \frac{k_o
      |\bx-\bx'|^2}{2  (z-z')} \Big)    .
\end{equation}
In a random medium, the 
first two moments of the random fundamental solution have the following expressions.

\begin{prop}
\label{prop:parax2}%
The first order-moment of the random  fundamental solution exhibits  damping
(for $z > z'$):
\ban
\EE \big[ \hat{g}\big( (\bx,z),(\bx',z') \big)   \big] 
&=&
\hat{g}_0\big( (\bx,z),(\bx',z') \big) 
 \exp \Big( -\frac{\gamma_0({\bf 0}) k_o^2 (z-z')}{8}  \Big)  ,
  \label{eq:mom1parax1}
\ean
where $\gamma_0$ is given by (\ref{def:gamma0}).

The second order-moment of the random  fundamental solution exhibits spatial decorrelation:
\ba
\nonumber 
&&  \hspace*{-2cm} \EE \big[ \hat{g}\big( (\bx_1,z),(\bx',z') \big) 
\overline{\hat{g}\big( (\bx_2,z),(\bx',z') \big)} \big]  =
\hat{g}_0\big( (\bx_1,z),(\bx',z') \big) 
\overline{\hat{g}_0\big( (\bx_2,z),(\bx',z') \big)} \\
 &&  \times
 \exp \Big( -  \frac{ \gamma_2(\bx_1-\bx_2) k_o^2 (z-z')}{4}   \Big)   , 
 \label{eq:mom2parax1}
\ea
where 
\begin{equation}
\gamma_2(\bx)= \int_0^1 \gamma_0({\bf 0}) -\gamma_0(\bx s) \d s .
\end{equation}
\end{prop}

These are classical results  (see \cite[Chapter 20]{ishimaru} and \cite{garniers2})
once the paraxial and white-noise approximations have been proved to be correct, as is the case here.
The result on the first-order moment shows that any coherent wave imaging method 
based on the mean field 
cannot give good images if the propagation distance is larger than the scattering mean free path  
\begin{equation}
\label{def:lsca:parax}
\ell_{\rm sca} = \frac{8 }{\gamma_0({\bf 0}) k_o^2},
\end{equation}
because the coherent wave components are then exponentially damped.
This is the situation we have in mind in this paper.

However, here the key quantity of interest is the intensity covariance function,
which means that we need to understand the behavior of the fourth-order moment of the field. We explain this next. 

\subsection{The statistical intensity covariance function}
In our paper the quantities of interest are the mean intensity 
\begin{equation}
{\cal I}_\br(\bx_0) =  \EE \big[ |E_{\br}(\bx_0)|^2\big]
\end{equation}
and the statistical intensity covariance function
\begin{equation}
{\cal C}_{\br,\br'} (\bx_0,\bx_0') = 
\EE \big[ |E_{\br} (\bx_0)|^2 |E_{\br'}(\bx_0')|^2 \big] - 
\EE \big[ |E_{\br}(\bx_0)|^2\big]\EE\big[ |E_{\br'}(\bx_0')|^2 \big]  .
\end{equation}
We remark that the statistical intensity covariance function is general in that   
 we have two,  in general different,  observation points  $\bx_0$ and $\bx_0'$,
 while in the  kernel in \eqref{def:intcor} the quadratic intensity term 
is evaluated at a common observation $\bx_0$.
We will discuss below,  in Section \ref{sec:Crho},  the measured intensity covariance function 
introduced in  \eqref{def:Crrrho0} and how it relates to the  mean intensity
and the  statistical intensity covariance function that we discuss here.

\begin{prop}
The second moment of the intensity can be expressed   as
\ba
\nonumber
\EE \big[ |E_{{\bf 0}}(\bx_0)|^2 |E_{\br}(\bx_0')|^2 \big] 
&=&\frac{1}{(2\pi)^{4d}} 
\iint_{\RR^{4d}} 
e^{i \bzeta_1 \cdot(\bx_0+\bx_0')+ i \bzeta_2 \cdot(\bx_0-\bx_0')}\\
&\times& 
 \hat{\mu}_\br (  \bxi_1,\bxi_2, \bzeta_1, \bzeta_2,{\ell})  
\d \bxi_1\d \bxi_2 \d \bzeta_1 \d  \bzeta_2  , \label{eq:2M}
\ea
where   $\hat{\mu}_\br$  satisfies
\begin{eqnarray}
\nonumber
&& 
\frac{\partial \hat{\mu}_\br}{\partial z} + \frac{i}{k_o} \big( \bxi_1\cdot \bzeta_1+   \bxi_2\cdot \bzeta_2\big) \hat{\mu}_\br
=
\frac{k_o^2}{4 (2\pi)^d} 
\int_{\RR^d} \hat{\gamma}_0(\bk) \Big[  
 \hat{\mu}_\br (  \bxi_1-\bk, \bxi_2-\bk, \bzeta_1, \bzeta_2)  \\
\nonumber
&& \quad  + 
 \hat{\mu}_\br (  \bxi_1-\bk,\bxi_2,  \bzeta_1, \bzeta_2-\bk)    
 +
 \hat{\mu}_\br (  \bxi_1+\bk, \bxi_2-\bk, \bzeta_1, \bzeta_2)   \\
\nonumber
&&  \quad 
+ 
 \hat{\mu}_\br (  \bxi_1+\bk,\bxi_2, \bzeta_1,  \bzeta_2-\bk)    -
2 \hat{\mu}_\br(\bxi_1,\bxi_2, \bzeta_1, \bzeta_2) \\
&&  \quad 
-
 \hat{\mu}_\br (  \bxi_1,\bxi_2-\bk, \bzeta_1, \bzeta_2-\bk)  
- \hat{\mu}_\br (  \bxi_1,\bxi_2+\bk,  \bzeta_1, \bzeta_2-\bk) 
\Big] \d \bk ,
\label{eq:fouriermom0}
\end{eqnarray}
starting from 
\begin{eqnarray}
\nonumber
&& \hat{\mu}_\br (  \bxi_1,\bxi_2, \bzeta_1, \bzeta_2,z=0)  
=
\hat{U}\Big( \frac{\bxi_1+\bxi_2+\bzeta_1+\bzeta_2}{2}\Big)
\overline{\hat{U}}\Big( \frac{\bxi_1+\bxi_2-\bzeta_1-\bzeta_2}{2}\Big) \\
&& \hspace*{-0.15in} \times 
\hat{U}\Big( \frac{\bxi_1-\bxi_2+\bzeta_1-\bzeta_2}{2}\Big)
\overline{\hat{U}}\Big( \frac{\bxi_1-\bxi_2-\bzeta_1+\bzeta_2}{2}\Big) \exp \big( i \br \cdot (\bzeta_2-\bzeta_1)\big).
\end{eqnarray}
 \end{prop}
No closed-form expression of the fourth moment of the field or of the second moment of the intensity is available,
but it is possible to get explicit expressions in two asymptotic regimes, the scintillation regime and the spot-dancing regime,
which correspond to the cases where the correlation radius of the medium is smaller (resp. larger) than the incident field radius and we discuss these in the next two sections. 
As we show in Section \ref{sec:Crho} the speckle imaging scheme considered here will
work well in the scintillation regime, however, as follows from the discussion in Section
\ref{sec:spot} not well in the spot dancing regime.  

\debproof
Note first that, by statistical transverse stationarity of the random  medium,
we have
\begin{equation}
{\cal C}_{\br,\br'}(\bx_0,\bx_0')  = {\cal C}_{{\bf 0},\br'-\br} \big(\bx_0-\br,\bx_0'-\br \big)  ,
\end{equation}
 It is therefore sufficient to study ${\cal C}_{{\bf 0},\br} \big(\bx_0 ,\bx_0'  \big)$.
We can write
\begin{eqnarray}
\EE \big[ |E_{{\bf 0}}(\bx_0)|^2 |E_{\br}(\bx_0')|^2 \big] =
{\cal M}_{\br}( \bx_0,\bx_0',\bx_0,\bx_0',{\ell}) ,
\end{eqnarray}
where 
we find using \eqref{eq:IS} and    the
It\^o theory for Hilbert-space valued random processes \cite{kunita}  that 
the fourth-order moment 
${\cal M}_{\br}( \bx_1,\bx_2,\by_1,\by_2,z)$ 
is solution of
\begin{eqnarray}\label{eq:M}
&&\frac{\partial {\cal M}_\br}{\partial z} = \frac{i}{2k_o}  \Big(  \Delta_{\bx_1}+\Delta_{\bx_2}
-  \Delta_{\by_1} -\Delta_{\by_2}\Big) {\cal M}_\br+ \frac{k_o^2}{4} {\cal U} \big( \bx_1,\bx_2, \by_1,\by_2\big)
{\cal M}_\br  , \\
&& {\cal M}_{\br}( \bx_1,\bx_2,\by_1,\by_2,z=0) 
=
U (\bx_1) \overline{U (\by_1)}
U_{\br}(\bx_2) \overline{U_{\br}(\by_2)},
 \end{eqnarray}
with the generalized potential
\begin{eqnarray}
&&   \hspace*{-.9cm}
{\cal U}\big( \bx_1,\bx_2, \by_1,\by_2 \big)  
=
\sum_{j,l=1}^2 \gamma_0(\bx_j-\by_l) 
- \gamma_0( \bx_1-\bx_2)
-  \gamma_0( \by_1-\by_2) -
2\gamma_0({\bf 0})   ,
\end{eqnarray}
and where $U$ is the shape of the mask as in  Eq. (\ref{eq:inc}).

We parameterize  the four points 
$\bx_1,\bx_2,\by_1,\by_2$  in the special way:
\begin{eqnarray}
\label{eq:reliexr1}
\bx_1 = \frac{\br_1+\br_2+\bq_1+\bq_2}{2}, \quad \quad 
\by_1 = \frac{\br_1+\br_2-\bq_1-\bq_2}{2}, \\
\bx_2 = \frac{\br_1-\br_2+\bq_1-\bq_2}{2}, \quad \quad 
\by_2 = \frac{\br_1-\br_2-\bq_1+\bq_2}{2}.
\label{eq:reliexr2}
\end{eqnarray}
We denote by $\mu_\br$ the fourth-order moment in these new variables:
\begin{equation}
\mu_\br (\bq_1,\bq_2,\br_1,\br_2,z) := 
{\cal M}_\br  (
\bx_1  ,
\bx_2 ,
\by_1 ,  
\by_2 ,z
) ,
\end{equation}
with $\bx_1,\bx_2,\by_1,\by_2$ given by (\ref{eq:reliexr1}-\ref{eq:reliexr2}) in terms of $\bq_1,\bq_2,\br_1,\br_2$.
The Fourier transform (in $\bq_1$, $\bq_2$, $\br_1$, and $\br_2$) of the fourth-order moment
is defined by:
\ba
\nonumber
\hat{\mu}_\br(\bxi_1,\bxi_2,\bzeta_1,\bzeta_2,z) 
&=& 
\iint_{\RR^{4d}} {\mu}_\br(\bq_1,\bq_2,\br_1,\br_2,z)  \\
&& \hspace*{-1.3in}
\times
\exp  \big(- i\bq_1 \cdot \bxi_1- i\bq_2 \cdot \bxi_2- i\br_1\cdot \bzeta_1- i\br_2\cdot \bzeta_2\big) \d \bq_1\d \bq_2 
\d \br_1 \d \br_2 \label{eq:fourier} 
.
\ea
When then arrive at Eq. (\ref{eq:2M}) using Eq. (\ref{eq:M}) and the Fourier transform. 
\finproof

\section{The scintillation regime}
\label{sec:scin}%
The scintillation regime is a  physically important 
regime corresponding to order one
relative fluctuations for the intensity.
The scintillation regime is valid if {the white-noise paraxial regime (Definition \ref{def:par}) is valid,
and, additionally,  the correlation radius of the medium fluctuations
(that determines the transverse correlation radius of the Brownian field in the It\^o-Schr\"odinger equation)
is smaller than the incident field radius.
The standard deviation of the Brownian field then needs to be relatively small and the propagation distance needs to be relatively large
to observe an effect of order one.}
More precisely,  we define the scintillation regime as follows.
\begin{definition}
\label{def:scint}
Consider the paraxial regime of Definition  \ref{def:par} 
so that the evolution of the field amplitude is governed
by the It\^o-Schr\"odinger equation (\ref{eq:IS}).
In the scintillation  regime,
 \begin{enumerate}
 \item
the covariance function $\gamma_0^\eps$  has an amplitude of order $\eps$:
\begin{equation}
\label{sca:sci}
\gamma_0^\eps(\bx)= \eps \gamma_0 (\bx) , 
\end{equation}
\item
the radius of the incident field and the vector shift are of order $1/\eps$:  
\begin{equation}
\label{def:feps}
U^\eps_\br(\bx) = U\big( \eps( \bx -\br)\big) ,  
\end{equation}
\item
the propagation distance is of order of $1/\eps$:
\begin{equation}
\label{def:Leps}
\ell^\eps =\frac{L}{\eps} ,
\end{equation}
\end{enumerate}
for a small dimensionless $\eps$.
 \end{definition}
Note that this problem was analyzed in \cite{garniers4} when $\br={\bf 0}$ and 
$U$ has a Gaussian profile. The following proposition \ref{prop:sci1} is an extension of this original result.

\subsection{The fourth-order moment of the transmitted field}
Let us denote the rescaled function
\begin{equation}
\label{eq:renormhatM2}
\tilde{\mu}_\br^\eps (\bxi_1,\bxi_2,\bzeta_1,\bzeta_2,z) := 
\hat{\mu}_\br \Big(\bxi_1,\bxi_2,\bzeta_1,\bzeta_2 , \frac{z}{\eps} \Big)
 \exp \Big( \frac{i  z}{k_o \eps} (\bxi_2 \cdot \bzeta_2  +   \bxi_1 \cdot \bzeta_1) \Big) .
\end{equation}
Our goal is to study the asymptotic behavior of $\tilde{\mu}_\br^\eps$ as $\eps \to 0$.
We have the following result, which shows that $\tilde{\mu}_\br^\eps$ exhibits a multi-scale behavior
as $\eps \to 0$, with some components evolving at the scale $\eps$ and 
some components evolving at the order one scale. 
The proof is similar to  the one of Proposition 1 in  \cite{garniers4}.
In \cite{garniers4} we used a Gaussian source profile  
while we here need to extend the result to the case of a general incident 
field, thus the calculus of the Gaussian for the source shapes do not apply 
directly as before.  However, the main steps of the proof remain unchanged and 
we obtain the following proposition.
\begin{prop}
\label{prop:sci1}%
In the scintillation regime of Definition \ref{def:scint}, 
if $\gamma_0\in L^1(\RR^d)$ and $\gamma_0({\bf 0})<\infty$, then
the function $\tilde{\mu}^\eps_\br(\bxi_1,\bxi_2, \bzeta_1,\bzeta_2,z ) $ 
can be expanded as
\begin{eqnarray}
\nonumber
&&
 \tilde{\mu}_\br^\eps(\bxi_1,\bxi_2,  \bzeta_1,\bzeta_2,z )  =
\frac{K(z)}{\eps^{4d}}
\hat{U}  \Big( \frac{\bxi_1+\bxi_2+\bzeta_1+\bzeta_2}{2 \eps}\Big)
\overline{\hat{U}} \Big( \frac{\bxi_1+\bxi_2-\bzeta_1-\bzeta_2}{2 \eps}\Big)\\
\nonumber
&& \hspace*{0.4in}\times
\hat{U}  \Big( \frac{\bxi_1-\bxi_2+\bzeta_1-\bzeta_2}{2 \eps}\Big)
\overline{\hat{U}} \Big( \frac{\bxi_1-\bxi_2-\bzeta_1+\bzeta_2}{2 \eps}\Big)
\exp \Big( i \br \cdot \frac{\bzeta_2-\bzeta_1}{\eps}\Big)
 \\
\nonumber
&& 
\quad
+
\frac{K(z)}{\eps^{3d}} 
\hat{V}_{\bf 0} \Big(\frac{\bzeta_2+\bzeta_1}{\eps}\Big)
\hat{U}  \Big( \frac{\bxi_1-\bxi_2+\bzeta_1-\bzeta_2}{2 \eps}\Big)
\overline{\hat{U}} \Big( \frac{\bxi_1-\bxi_2-\bzeta_1+\bzeta_2}{2 \eps}\Big)\\
\nonumber
&& \hspace*{0.4in}\times
\exp \Big( i \br \cdot \frac{\bzeta_2-\bzeta_1}{\eps}\Big)
A\big(\frac{\bxi_2+\bxi_1}{2} ,\frac{\bzeta_2 + \bzeta_1}{\eps} ,z \big) \\
\nonumber
&&  
\quad
+
\frac{K(z)}{\eps^{3d}} 
\overline{\hat{V}_{\bf 0}} \Big(\frac{\bzeta_2-\bzeta_1}{\eps}\Big)
\hat{U}  \Big( \frac{\bxi_1+\bxi_2+\bzeta_1+\bzeta_2}{2 \eps}\Big)
\overline{\hat{U}} \Big( \frac{\bxi_1+\bxi_2-\bzeta_1-\bzeta_2}{2 \eps}\Big)\\
\nonumber
&& \hspace*{0.4in}\times
\exp \Big( i \br \cdot \frac{\bzeta_2-\bzeta_1}{\eps}\Big)
A \big(\frac{\bxi_2-\bxi_1}{2} ,\frac{\bzeta_2- \bzeta_1}{\eps} ,z \big) \\
\nonumber
&&  
\quad
+
\frac{K(z)}{\eps^{3d}} 
\hat{V}_{\br} \Big(\frac{\bxi_2+\bzeta_1}{\eps}\Big)
\hat{U}  \Big( \frac{\bxi_1-\bzeta_2+\bzeta_1-\bxi_2}{2 \eps}\Big)
\overline{\hat{U}} \Big( \frac{\bxi_1-\bzeta_2-\bzeta_1+\bxi_2}{2 \eps}\Big)\\
\nonumber
&& \hspace*{0.4in}\times
\exp\Big( i \frac{\br}{2} \cdot \frac{\bzeta_2-\bxi_1-2\bzeta_1}{\eps}\Big)
A\big( \frac{\bzeta_2+\bxi_1}{2} ,\frac{\bxi_2+ \bzeta_1}{\eps} ,z \big) \\
\nonumber
&& 
\quad
+
\frac{K(z)}{\eps^{3d}} 
\overline{\hat{V}_{\br}} \Big(\frac{\bxi_2-\bzeta_1}{\eps}\Big)
\hat{U}  \Big( \frac{\bxi_1+\bxi_2+\bzeta_1+\bzeta_2}{2 \eps}\Big)
\overline{\hat{U}} \Big( \frac{\bxi_1-\bxi_2-\bzeta_1+\bzeta_2}{2 \eps}\Big)\\
\nonumber
&& \hspace*{0.4in}\times
\exp\Big( i \frac{\br}{2} \cdot \frac{\bzeta_2+\bxi_1-2\bzeta_1}{\eps}\Big)
A\big(  \frac{\bzeta_2-\bxi_1}{2} ,\frac{\bxi_2- \bzeta_1}{\eps}  ,z \big) \\
\nonumber
&& 
\quad
+\frac{K(z)}{\eps^{2d}} 
\hat{V}_{\bf 0} \Big( \frac{\bzeta_2+\bzeta_1}{\eps}\Big)
\overline{\hat{V}_{\bf 0} }\Big( \frac{\bzeta_2-\bzeta_1}{\eps}\Big) \\
\nonumber
&& \hspace*{0.4in}\times
\exp\Big( i \br \cdot \frac{\bzeta_2-\bzeta_1}{\eps}\Big)
A \big( \frac{\bxi_2+\bxi_1}{2},   \frac{\bzeta_2+ \bzeta_1}{\eps} ,z \big)
A \big( \frac{\bxi_2-\bxi_1}{2},   \frac{\bzeta_2- \bzeta_1}{\eps} ,z  \big) \\
\nonumber
&&
\quad
 + \frac{K(z)}{\eps^{2d}} 
 \hat{V}_{\br} \Big( \frac{\bxi_2+\bzeta_1}{\eps}\Big)
\overline{ \hat{V}_\br } \Big( \frac{\bxi_2-\bzeta_1}{\eps}\Big) \\
\nonumber
&& \hspace*{0.4in}\times
\exp\Big(- i \br \cdot \frac{\bzeta_1}{\eps}\Big)
A \big( \frac{\bzeta_2+\bxi_1}{2},  \frac{\bxi_2+ \bzeta_1}{\eps} ,z  \big)
A \big( \frac{\bzeta_2-\bxi_1}{2},  \frac{\bxi_2- \bzeta_1}{\eps} ,z \big)
\\
&& 
\quad
 + R^\eps  ( \bxi_1,\bxi_2 ,  \bzeta_1 ,\bzeta_2 ,z )   ,
\label{eq:propsci11}
\end{eqnarray}
where the functions $K$ and $A$  are defined by
\begin{eqnarray}
\label{def:K}
K(z) &:=&  
 \exp\Big(- \frac{k_o^2}{2} \gamma_0({\bf 0}) z\Big) , \\
\nonumber
A(\bxi,\bzeta,z)  &:=& 
 \frac{1}{(2\pi)^d}
 \int_{\RR^d}  \Big[  \exp \Big( \frac{k_o^2}{4} \int_0^z \gamma_0 \big( \bx + \frac{ \bzeta}{k_o} z' \big) \d z' \Big) -1\Big]\\
 && \times
   \exp \big( -i \bxi\cdot \bx  \big)
 \d \bx  ,   
\label{def:A}
\end{eqnarray}
the function $\hat{V}_\br$ is
\begin{equation}
\hat{V}_\br(\bzeta) = \int \hat{U} \big(\bk + \frac{\bzeta}{2}\big)
\overline{ \hat{U} }  \big(\bk - \frac{\bzeta}{2}\big) \exp\big( i \bk \cdot \br \big)
\d \bk,
\end{equation}
and the function $R^\eps $ satisfies
\ban
\sup_{z \in [0,{\ell}]} \| R^\eps (\cdot,\cdot,\cdot,\cdot, z ) \|_{L^1(\RR^d\times \RR^d\times \RR^d\times \RR^d)} 
\stackrel{\eps \to 0}{\longrightarrow}  0  .
\ean
\end{prop}
It is shown in \cite{garniers4} that the function $\bxi \to A(\bxi,\bzeta,z)$
 belongs to $L^1(\RR^d)$ 
and  its $L^1$-norm $\| A(\cdot,\bzeta,z)\|_{L^1(\RR^d)}$ is bounded uniformly in 
$\bzeta \in \RR^d$ and $z\in [0,{\ell}]$. 
It follows that   all terms in the expansion (except the remainder $R^\eps$) 
have $L^1$-norms of  order one
when $\eps \to 0$.

\subsection{The statistical intensity covariance function}
We will here characterize  the mean intensity and the statistical intensity covariance function  when
\begin{equation}
\label{eq:paramx0} 
\bx_0= \frac{\bX_0}{\eps}+\frac{\bY_0}{2},\quad \quad 
\bx_0'= \frac{\bX_0}{\eps}-\frac{\bY_0}{2} ,  \quad \quad 
\br= \frac{\bR}{\eps},
\quad \quad \br'= \frac{\bR'}{\eps}.
\end{equation}
Here  coordinates in capital letters are of  order one with respect to $\eps$,
so that $\bX_0, \bR, \bR'$  are rescaled lateral  coordinates and
$\bY_0$ is the observation
offset in original coordinates. 
{This means that we consider the intensity covariance function for mid-points located within the beam whose radius is large 
(of order $\eps^{-1}$) but we look at offsets that are small (of the order of the correlation length of the medium, that is of order one).
The motivation for this parameterization is indeed 
that the intensity distribution decorrelates for the  observation offset  on this scale, while 
we will see that the intensity covariance function as a function of $\br$ and $\br'$ varies naturally at the scale $\eps^{-1}$.}
Recall that, by (\ref{def:Leps}),   $L$ is the propagation distance from the mask to the camera
in rescaled longitudinal  coordinates.
We have the following result.

\begin{prop}\label{prop:2}
In the scintillation regime,
we  have in the limit  $\eps \to 0$
 \ba
\nonumber
{\cal I}_\br (\bx_0) &=&
\frac{1}{(2\pi)^d} \int_{\RR^d}
 \Big( \int_{\RR^d} |U(\bX-\bR)|^2 \exp \big( - i\bzeta\cdot \bX \big) \d \bX\Big) \\
&  \times &
\exp \big( i\bzeta\cdot \bX_0 \big)
 \exp \Big( \frac{k_o^2}{4} \int_0^L \gamma_0 \big( \frac{ \bzeta}{k_o} z\big) - \gamma_0({\bf 0})  \d z \Big)
 \d \bzeta   ,
 \label{eq:m}
\ea
and 
\ba\nonumber
 && \hspace*{-2cm}  {\cal C}_{\br,\br'}(\bx_0,\bx_0')  = 
\Big|
\frac{1}{(2\pi)^d} \int_{\RR^d}
 \Big( \int_{\RR^d} U\big(\bX+\frac{\bR'-\bR}{2}\big)\overline{U}\big(\bX-\frac{\bR'-\bR}{2}\big)
  \exp \big( - i\bzeta\cdot \bX \big)\d \bX\Big) \\
\nonumber
& &  \times
\exp \Big( i\bzeta\cdot \big(\bX_0- \frac{\bR+\bR'}{2}\big)\Big)
 \exp \Big( \frac{k_o^2}{4} \int_0^L \gamma_0 \big( \frac{ \bzeta}{k_o} z-\bY_0 \big) -\gamma_0({\bf 0})  \d z \Big)
 \d \bzeta
\Big|^2 \\
\nonumber
&  &  \hbox{}  -\Big|
\frac{1}{(2\pi)^d} \int_{\RR^d}
 \Big( \int_{\RR^d} U\big(\bX+\frac{\bR'-\bR}{2}\big)
 \overline{U}\big(\bX-\frac{\bR'-\bR}{2}\big) \exp \big( - i\bzeta\cdot \bX \big)\d \bX\Big) \\  \label{eq:c}
&  & \times
\exp \Big( i\bzeta\cdot \big(\bX_0- \frac{\bR+\bR'}{2} \big)\Big)
\exp\Big(- \frac{k_o^2}{4} \gamma_0({\bf 0}) L \Big) 
 \d \bzeta
\Big|^2   .
\ea
\end{prop}
\debproof
   The result  follows from Proposition \ref{prop:sci1}.      
\finproof
 
In order to get explicit expressions for the quantity of
interest it is convenient to introduce the  
 strongly scattering regime defined as follows.
 Recall that the scattering mean free path $\ell_{\rm sca}$ is defined by (\ref{def:lsca:parax})).
\begin{definition}
\label{def:sto}
In the strongly scattering regime, we have $L/\ell_{\rm sca} \gg 1$
and {the fluctuations of the random medium are smooth so that} the function $\gamma_0$ can be expanded as
\begin{equation}
\label{eq:expandgamma0}
\gamma_0(\bx) = \gamma_0({\bf 0}) - \frac{1}{2} \bar{\gamma}_2 |\bx|^2 + o(|\bx|^2),
\end{equation}
for $\bx$ {smaller than the correlation length of the medium (i.e., the width of $\gamma_0$)}.
 \end{definition}
 This corresponds to large,  but smooth,  medium fluctuations. 
 We can now identify simplified expressions  for the mean  intensity and the
 intensity  covariance function. 
 
\begin{lemma}\label{lem:1}
Assume the scintillation and   strongly scattering regime, then   
we have
\ba
{\cal I}_\br (\bx_0)  =
 \frac{6^{d/2}}{(\pi \bar{\gamma}_2 L^3)^{d/2}}  \int_{\RR^d}
 |U(\bX)|^2  \exp \Big( -  \frac{6 | \bX- \bX_0+\bR |^2}{\bar{\gamma}_2 L^3}  \Big)
 \d \bX
 \label{eq:Irx0}
\ea
and
 \begin{eqnarray}
\nonumber
&&{\cal C}_{\br,\br'}(\bx_0,\bx_0') = 
 \frac{6^d}{(\pi \bar{\gamma}_2 L^3)^d} 
\Big|
   \int_{\RR^d} U\big(\bX+\frac{\bR'-\bR}{2}\big)\overline{U}\big(\bX-\frac{\bR'-\bR}{2}\big)
  \\
  \nonumber
&&   ~~  \times
\exp \Big( -  \frac{6 | \bX- \bX_0+ \frac{\bR+\bR'}{2} |^2}{\bar{\gamma}_2 L^3} 
-i \frac{3 k_o}{2L} \bY_0 \cdot \big( \bX- \bX_0 + \frac{\bR+\bR'}{2}  \big) \Big)
 \d \bX
\Big|^2 \\
&&  ~~  \times\exp\Big( -\frac{\bar{\gamma}_2 k_o^2 L }{16} |\bY_0|^2 \Big).
\label{eq:covtot0}
\end{eqnarray}
\end{lemma} 
 
\debproof 
{In the scintillation regime, we can write by Eq.~(\ref{eq:expandgamma0}):
$$
\exp\Big( - \big(1-\frac{\gamma_0(\bx)}{\gamma_0({\bf 0})}\big) \frac{L}{\ell_{\rm sca}} \Big) \simeq
\exp
 \Big( - \frac{ \bar{\gamma}_2 L}{2\gamma_0({\bf 0})\ell_{\rm sca}} |\bx|^2\Big) ,
$$
since this is true for $|\bx|$ smaller than the correlation length,
moreover,  since this is also true for $|\bx|$ of the order of or larger than the correlation length
in the sense that the two members of the equations are exponentially small in $L/\ell_{\rm sca}$.
}
It then follows from Proposition \ref{prop:2} that   
  the mean intensity is
 \ba
\nonumber
{\cal I}_\br (\bx_0) &=&
\frac{1}{(2\pi)^d} \int_{\RR^d}
 \Big( \int_{\RR^d} |U(\bX-\bR)|^2 \exp \big( - i\bzeta\cdot \bX \big) \d \bX\Big) \\
& &   \times
\exp \big( i\bzeta\cdot \bX_0 \big)
 \exp \Big( -  \frac{\bar{\gamma}_2 L^3}{24} |\bzeta|^2 \Big)
 \d \bzeta
\ea
 and 
 the intensity covariance function is
  \begin{eqnarray*}
\nonumber
&& 
{\cal C}_{\br,\br'}(\bx_0,\bx_0') =  \\ &&  \hspace{-1.7cm}
\Big|
\frac{1}{(2\pi)^d} \int_{\RR^d}
 \Big( \int_{\RR^d} U\big(\bX+\frac{\bR'-\bR}{2}\big)\overline{U}\big(\bX-\frac{\bR'-\bR}{2}\big)
  \exp \big( - i\bzeta\cdot \bX \big)\d \bX\Big) \\ \nonumber 
&&  \hspace{-1.7cm}   \times
\exp \Big( i\bzeta\cdot \big(\bX_0- \frac{\bR+\bR'}{2} \big)\Big)
\exp \Big( -\frac{\bar{\gamma}_2 L^3}{24} |\bzeta|^2 + \frac{\bar{\gamma}_2 k_oL^2}{8} \bzeta \cdot \bY_0 
-\frac{\bar{\gamma}_2 k_o^2 L }{8} |\bY_0|^2 \Big)
 \d \bzeta
\Big|^2 . \nonumber
\end{eqnarray*}
The lemma then follows after integrating in $\bzeta$. 
\finproof

The beam radius enhancement  due to scattering 
in a random medium with thickness $L$ is given by \cite[Eq.~(74)]{garniers4}:
\ba\label{eq:cA}
\cA_L  := \bar{\gamma}_2^{1/2} L^{3/2}  /(\eps \sqrt{6})   .
\ea
Let us assume a regime of  large enhanced aperture defined as follows.
\begin{definition}\label{def:enh}
In the large enhanced aperture regime,
 the radius of the  incident field $U$,   
the radius and center point of the camera,  and the 
shifts $|\br|,|\br'|$ are small  relative to the  beam radius 
enhancement $\cA_L$. 
\end{definition}
 {As we show below this is the configuration in which 
 the intensity covariance function has a simple form and 
 the profile of the 
incident field can be explicitly extracted, because one can extract a large range of values in $\br-\br'$ 
of the intensity covariance function}.  
We can also address the general situation, 
albeit with less  explicit expressions,  
and we do so  in Remark \ref{remark2}. 
It follows from Lemma  \ref{lem:1} that   in the  large enhanced aperture
regime we have the following result.
\begin{lemma}
In the scintillation,    strongly scattering,  and large enhanced aperture regime, 
 the mean intensity is constant over the camera:
\begin{equation}
{\cal I}_\br (\bx_0) =  \frac{6^{d/2}}{(\pi \bar{\gamma}_2 L^3)^{d/2}}  \int_{\RR^d}
 |U(\bX)|^2   \d \bX  ,
  \label{eq:Irx0b}
\end{equation}
and the intensity covariance function is
 \ba
\nonumber
{\cal C}_{\br,\br'}(\bx_0,\bx_0') &=&
 \frac{6^d}{(\pi \bar{\gamma}_2 L^3)^d} 
\Big|
   \int_{\RR^d} U\big(\bX+\frac{\bR'-\bR}{2}\big)\overline{U}\big(\bX-\frac{\bR'-\bR}{2}\big)
 \d \bX
\Big|^2 \\
& &  
\times \exp \Big(  -\frac{\bar{\gamma}_2 k_o^2 L }{16} |\bY_0|^2 \Big).
\label{eq:Crrp2}
\ea
\end{lemma}
Thus, the intensity covariance function does not depend on the mid observation point $\bX_0$ and on the shift mid-point $(\bR+\bR')/2$, 
but it decays as a function of the shift offset 
$\bR'-\bR$ on a scale length that is of the order of the incident field radius 
 in a way that makes it possible to reconstruct the incident field.
It was shown in  \cite[Proposition 6.3]{garniers3}  and  also in  \cite[Eq.~(75)]{garniers4} that 
\begin{equation}
\label{def:rhoL}
\rho_L := \frac{2}{ \sqrt{\bar{\gamma}_2 k_o^2  L}}     
\end{equation}
is the typical correlation radius of the speckle pattern generated by a {plane wave} going through a random medium with thickness $L$.
 We can see from (\ref{eq:Crrp2}) that the intensity covariance function 
   decays as a whole with the observation offset $\bY_0$ on a scale length equal to $\rho_L$. 
  Indeed,   the intensity covariance function
decays on the  scale $\rho_L$ with respect to observation offset because 
the speckle pattern, that is the intensity fluctuations,   decorrelates on this scale.

\subsection{Extraction of the incident field profile}
\label{sec:Crho}%
The empirical intensity covariance function is given by (\ref{def:Crrrho0}).
If  the radius $r_A=R_A/\eps$   of the camera is large enough 
(more precisely, if condition (\ref{cond:sa}) below holds true, {which means that the camera covers many speckle spots}),
then the empirical intensity covariance function  is self-averaging and equal to
\ba
{\cal C}^{\rho_o}_{\br,\br'} = \frac{1}{(4\pi)^{d/2}\rho_o^d}
\int  {\cal C}_{\br,\br'}(\bx_0,\bx_0') \exp \Big( - \frac{|\bY_0|^2}{4 \rho_o^2} \Big)  \d \bY_0  .
\label{eq:Crrpa}
\ea
{This is because $ \frac{1}{|A_o|}
\int_{A_o} \cdots \d \bx_0$ in (\ref{def:Crrrho0}) becomes equal to $\EE[\cdots]$ by the law of large numbers.}
 Here we used   the parameterization \eqref{eq:paramx0} and the expressions
\eqref{eq:m} and \eqref{eq:c} which show in particular that the mean intensity varies
on the  slow  scale $\eps^{-1}$ relative to the characteristic speckle size, the scale of 
decorrelation of the intensities. 
We also remark that the condition (\ref{cond:sa}) 
 means that the camera has many pixels and  also observes
many speckle spots.
The result \eqref{eq:Crrpa} is valid in the general
 scintillation case.  We next present the main result of the paper.
 
  \begin{prop}\label{prop:M1} 
Assume the scintillation, strongly  scattering, and large enhanced aperture   
regimes of Definitions \ref{def:scint}, \ref{def:sto}, and  \ref{def:enh} 
 respectively. Moreover, assume that the radius of
the camera satisfies
\begin{equation}
\label{cond:sa}
R_A  \gg \sqrt{\rho_o^2 +\rho_L^2} . 
\end{equation}
Then the intensity covariance function is self-averaging and given by
\begin{eqnarray}
\nonumber
{\cal C}^{\rho_o}_{\br,\br'}
&=&{\cal Z}^{\rho_o} \Big|
\int_{\RR^d}
U \big(\bX+\frac{\bR-\bR'}{2} \big)\overline{U} \big(\bX-\frac{\bR-\bR'}{2}\big)  \d\bX
 \Big|^2 \\
&=&{\cal Z}^{\rho_o} \Big| \frac{1}{(2\pi)^d} 
\int_{\RR^d}
|\hat{U}(\bzeta)|^2 \exp \big( i\bzeta \cdot \frac{\bR'-\bR}{2} \big)  \d\bzeta
 \Big|^2 ,  \label{eq:Crf}
\end{eqnarray}
with 
\begin{equation}
{\cal Z}^{\rho_o} =   \frac{6^d}{(\pi \bar{\gamma}_2 L^3)^d}  \frac{1}{\big( 1+  {\rho_o^2}/{\rho_L^2} \big)^{d/2}} .
\end{equation}
\end{prop} 
The  multiplicative factor ${\cal Z}^{\rho_o}$ is the 
one associated with the mean square intensity in \eqref{eq:I2}  below 
  when the pixel size $\rho_o$ is smaller than the speckle size $\rho_L$.
This is exactly the formula (\ref{eq:pred}) predicted in \cite{webb14}.

\debproof
 We have assumed:\\
1)
scattering is strong $L/\ell_{\rm sca} \gg 1$, leading to \eqref{eq:expandgamma0};\\
2) the radius $R_A$ of the camera $A_o$ is smaller 
than $\bar{\gamma}_2^{1/2} L^{3/2}$;\\
3) 
the shifts  $|\bR|,|\bR'|$ and camera center point magnitude 
are smaller than $\bar{\gamma}_2^{1/2} L^{3/2}$.\\
Then the intensity covariance function has the form \eqref{eq:Crrp2}.
 Substituting (\ref{eq:Crrp2}) into \eqref{eq:Crrpa}  gives
\begin{equation} \hspace*{-1.5cm} 
{\cal C}^{\rho_o}_{\br,\br'} =  \frac{6^d}{(\pi \bar{\gamma}_2 L^3)^d}  
\frac{1}{\big( 1+ {\rho_o^2}/{\rho_L^2} \big)^{d/2}}
\Big|
   \int_{\RR^d} U\big(\bX+\frac{\bR'-\bR}{2}\big)\overline{U}\big(\bX-\frac{\bR'-\bR}{2}\big)
 \d \bX \Big|^2 .
 \label{eq:Crrp3}
\end{equation}
{The self-averaging can be considered as efficient because the amplitude of the main peak of the intensity covariance function is larger than the fluctuations
of the background. Indeed the background is the square of the mean intensity (see (\ref{eq:Irx0b})): }
\ba\label{eq:I2}
{\cal I}^2 =  \frac{6^{d}}{(\pi \bar{\gamma}_2 L^3)^{d}} \Big( \int_{\RR^d}
 |U(\bX)|^2   \d \bX \Big)^2 ,
\ea 
and its fluctuations are  of the order of ${\cal I}^2/ \sqrt{M}$ where $M$ is the number of speckle spots over which the 
averaging has been carried out, that is, $M=(R_A  / \rho_L)^d$.
 Note   that in the strongly scattering 
scintillation regime the field $E_\br(\bx)$ will,
from the point of view of the fourth moment, behave as a complex-valued circularly symmetric
Gaussian random variable \cite{garniers4}, which means in particular that
$\EE[ | E_\br(\bx) |^2 ]^2  =  {\rm Var} [ | E_\br(\bx) |^2 ]$. 
 
The amplitude of the main peak of the intensity covariance function is (by (\ref{eq:Crrp3}))
\ban
{\cal C}_{{\bf 0},{\bf 0}}^{\rho_o} =  \frac{1}{\big( 1+ {\rho_o^2}/{\rho_L^2} \big)^{d/2}} {\cal I}^2 .
\ean
The main peak can be clearly estimated if
$
{1}/(  1+ {\rho_o^2}/{\rho_L^2} )^{d/2}  \gg  {\rho_L^d}/{R_A^d}
$,
which leads to the condition (\ref{cond:sa}). 
\finproofo

\begin{remark}
The results presented in this section show that the intensity covariance function over incident field position
makes it possible to reconstruct the incident field. One may ask whether the intensity covariance function over 
transmitted field position may also possess this property. 
To answer this question we inspect
${\cal C}_{{\bf 0},{\bf 0}}(\bx_0,\bx_0')$.
As shown by (\ref{eq:Crrp2}), this intensity covariance over transmitted field position   in the strongly scattering regime is
\begin{equation}
{\cal C}_{{\bf 0},{\bf 0}}(\bx_0,\bx_0') =
 \frac{6^d}{(\pi \bar{\gamma}_2 L^3)^d} 
\Big|   \int_{\RR^d} |U(\bX)|^2  \d \bX
\Big|^2 
\exp \Big(  -\frac{ |\bY_0|^2}{4\rho_L^2} \Big) ,
\end{equation}
when $\bx_0$ and $\bx_o'$ are as in (\ref{eq:paramx0}).
There is, therefore, no way to reconstruct the incident field given this function only.
\end{remark}

\begin{remark}
\label{remark2}
 To be complete, let us now address the case when the radius of the incident field is of the order of or even
larger than the enhanced aperture $\bar{\gamma}_2^{1/2} L^{3/2}$.
Then we also consider a camera $A_o$ with a radius
 larger than $\bar{\gamma}_2^{1/2} L^{3/2}$ and shifts $|\bR|,|\bR'|$ larger than $\bar{\gamma}_2^{1/2} L^{3/2}$. 
Under these circumstances Eq.~(\ref{eq:Irx0}) shows that the mean intensity gives a blurred version of the incident field profile,
in the form of a convolution of $|U|^2$ with a Gaussian kernel of  width 
of order $\bar{\gamma}_2^{1/2} L^{3/2}$.
 The intensity covariance function (\ref{eq:covtot0}) depends on the mid point $\bX_0$.
Let us first consider the situation when we integrate 
the actual intensity covariance function with respect to
mid point $\bX_0$,
which gives
\begin{eqnarray}
\nonumber
&&\int {\cal C}_{\br,\br'}(\bx_0,\bx_0')  \d\bX_0 \\
\nonumber
&&= 
\frac{1}{(2\pi)^d} 
\int_{\RR^d} \Big| \int_{\RR^d} U \big(\bX +\frac{\bR'-\bR}{2} \big)
\overline{U} \big(\bX - \frac{\bR'-\bR}{2} \big) e^{-i \bzeta \cdot \bX} \d\bX \Big|^2 \\
&& \quad \times
\exp \Big(-\frac{\bar{\gamma}_2 L^3}{12}\big| \bzeta  - \frac{3k_o}{2L} \bY_0\big|^2 \Big) 
\d\bzeta \exp \Big( - \frac{|\bY_0|^2}{4\rho_L^2} \Big).
\end{eqnarray}
 Thus, with an offset $\bY_0$ in the observation points there
is a damping of the information on the scale $\rho_L$ due to
  decorrelation of the speckle pattern. We also have a damping
 of the information at spatial scales of $U$ that are larger than   
$\bar{\gamma}_2^{1/2} L^{3/2}$. 
The empirical intensity covariance function is self-averaging and equal to
(\ref{eq:Crrpa}), so that the integrated (in $\bX_0$) version is equal to
\begin{eqnarray}
\nonumber
&&\int {\cal C}_{\br,\br'}^{\rho_o}(\bx_0,\bx_0')  \d\bX_0 \\
\nonumber
&& =  
\frac{3^{d/2}}{[\pi \bar{\gamma}_2L^3 (1+ \rho_o^2 /\rho_L^2 )]^{d/2}}
 \iint 
  U \big(\bX +\frac{\bR'-\bR}{2} \big)\overline{U} \big(\bX - \frac{\bR'-\bR}{2} \big) \\
  &&  \times 
  \overline{U} \big(\bX' +\frac{\bR'-\bR}{2} \big) {U} \big(\bX' - \frac{\bR'-\bR}{2} \big) 
  \exp \Big( - \frac{|\bX-\bX'|^2}{2R_L^2}  \Big) \d \bX \d\bX',
\end{eqnarray}
where
\begin{equation}
R_L^2 = \frac{\bar{\gamma}_2 L^3}{6}
 \frac{1+   \rho_o^2 /\rho_L^2}{1+4  \rho_o^2 /\rho_L^2}  ,
\end{equation}
which is between $\bar{\gamma}_2 L^3/{24}$ and $\bar{\gamma}_2 L^3/{6}$.
This shows that the intensity covariance function is proportional to
(\ref{eq:Crrp3}) when the radius of the incident field is smaller than $R_L$, but it 
becomes blurred by the Gaussian convolution with radius $R_L$
when the radius of the incident field is larger.
\end{remark}

{
\begin{remark}
In this paper we assume that the phase of $U$ is known,
so that a phase-retrieval algorithm can be used to extract $U$ from $|\hat{U}|$.
This is the case in the experimental setting described at the beginning of Section \ref{sec:intcor}, as 
we assume that 1) the illumination is a plane wave and 2) the object is a mask. 
If the plane wave is normally incident, then the phase is zero (or constant). 
If the plane is obliquely incident with a known angle, then the phase is also known.
If the illumination phase is unknown, then it should still be possible -in principle- to reconstruct the complex profile $U$
provided one has a sufficiently strong support constraint as demonstrated in \cite{fienup87},
but this is not obvious.
\end{remark}
}

{
\begin{remark}\label{remark:49}
In this paper we assume that the medium is statistically homogeneous between the 
plane of the mask $z=0$ and the plane of the camera $z=\ell$.
One could  consider a more general situation in which there are three regions,
namely a random medium sandwiched in between two homogeneous media.
This situation  will be addressed in a further work but we may anticipate the contributions of
interesting phenomena such as the shower curtain effect \cite{ishimaru}.
\end{remark}
}

\section{The spot-dancing regime}
\label{sec:spot}%
The  spot-dancing  regime is valid if {the white-noise paraxial regime (Definition \ref{def:par}) is valid,
and, additionally,  the correlation radius of the medium fluctuations
(that determines the transverse correlation radius of the Brownian field in the It\^o-Schr\"odinger equation)
is larger than the incident field radius.
The standard deviation of the Brownian field then needs to be relatively large so that one can see an effect of order one.}
More precisely, we define the spot-dancing regime as follows.
\begin{definition}
\label{def:spo}
{Consider the paraxial regime of Definition  \ref{def:par} 
so that the evolution of the field amplitude is governed
by the It\^o-Schr\"odinger equation (\ref{eq:IS}).}
In the spot-dancing regime,
the covariance function $\gamma_0^\epsilon$ is of the form:
\begin{equation}
\label{eq:spregime}
\gamma_0^\epsilon (\bx) = \epsilon^{-2} \gamma_0 (\epsilon \bx)  ,
\end{equation}
for a small dimensionless parameter $\epsilon$, and the function $\gamma_0$ is smooth and can be expanded as (\ref{eq:expandgamma0}).
\end{definition}

We want to study the asymptotic behavior of 
the moments of the field in this regime, which is called the 
spot-dancing regime for reasons
that will become clear from the following discussion.
\begin{prop}\label{pro:spo}
In the spot-dancing regime we have   the following asymptotic description for the 
transmitted field   in distribution 
\begin{eqnarray}
\nonumber
{E}_\br(\bx) &=&{E}^0_\br(\bx -\bX_{\ell}) 
\exp \Big( - i \frac{ k_o\sqrt{\bar{\gamma}_2} \bW_{\ell} }{2}
\cdot (\bx-\bX_{\ell}) \Big) \\
&&
\times \exp \Big( i \frac{k_o \bar{\gamma}_2 }{8}
\big( {\ell} |\bW_{\ell}|^2 - \int_0^{\ell} |\bW_{z}|^2 \d z \big) \Big)  ,
\end{eqnarray}
{where $\bW_z$ is a standard $d$-dimensional Brownian motion,}
\begin{equation}
{E}^0_\br(\bx) =
\Big( \frac{k_o}{2\pi {\ell}} \Big)^{d/2}
\int_{\RR^d} U(\by-\br) \exp \Big( i \frac{k_o}{2{\ell}} |\bx-\by|^2 \Big)\d \by
\label{eq:fieldfar1}
\end{equation}
is the field that is observed when the medium is homogeneous and
\begin{equation}
\bX_z = \frac{\sqrt{\bar{\gamma}_2}}{2} \Big(\int_0^z \bW_{z'} \d z'- z \bW_z \Big)
= -\frac{\sqrt{\bar{\gamma}_2}}{2}   \int_0^z z' \d \bW_{z'}  
\end{equation}
is the random center of the field, that is a $\RR^d$-valued Gaussian process with mean zero and covariance
\begin{equation}
\EE \big[ \bX_z \bX_{z'}^T \big] = \frac{\bar{\gamma}_2 (z\wedge z')^3}{12} {\bf I}  .
\end{equation}
\end{prop}
In particular the intensity of the transmitted field is
\begin{equation}
| {E}_\br (\bx)|^2  = | {E}^0_\br (\bx - \bX_{\ell})|^2 =  | {E}^0_{\bf 0} (\bx - \br-\bX_{\ell})|^2  .
\label{eq:meanintspot1}
\end{equation}
This representation justifies the name ``spot-dancing regime": the transmitted intensity has the same transverse
profile as in a homogeneous medium, but its center is randomly shifted by the Gaussian process 
$\bX_z$.
Note that in this case, there is no statistical averaging when one considers the empirical 
intensity covariance  function (\ref{def:Crrrho0}), which is the random quantity equal to 
\begin{eqnarray}
\nonumber &&   \hspace*{-1cm} 
 C^{\rho_o}_{\br,\br'} 
= \frac{1}{|A_o|}
\int_{A_o} | {E}^0_{\bf 0} (\bx_0 - \br-\bX_{\ell})|^2| {E}^0_{\bf 0} (\bx_0 - \br'-\bX_{\ell})|^2 \d \bx_0 \\
&& \hspace*{-1cm} 
- \Big( \frac{1}{|A_o|}
\int_{A_o} | {E}^0_{\bf 0} (\bx_0 - \br-\bX_{\ell})|^2  \d \bx_0 \Big)
\Big( \frac{1}{|A_o|}
\int_{A_o} | {E}^0_{\bf 0} (\bx_0 - \br'-\bX_{\ell})|^2 \d \bx_0 \Big)  .
\end{eqnarray}
If the radius of the camera is larger than the radius of the incident field, moreover,   large relative
to $\sqrt{\bar{\gamma}_2 {\ell}^3}$,  the typical spot dancing shift,  and the shift $\br$, 
then the intensity covariance  function gives the autocovariance of the unperturbed intensity profile:
\begin{eqnarray}
&& \hspace*{-0.6in}
C^{\rho_o}_{\br,\br'} =\frac{1}{|A_o|}
\int_{\RR^d} | {E}^0_{\bf 0} (\bx_0  )|^2| {E}^0_{\bf 0} (\bx_0 - \br'+\br)|^2 \d \bx_0
- \Big( \frac{1}{|A_o|}
\int_{\RR^d} | {E}^0_{\bf 0} (\bx_0  )|^2  \d \bx_0 \Big)^2
.
\label{eq:expresCrhospot2}
\end{eqnarray}
Therefore, in the spot-dancing regime, the random medium does not modify the intensity covariance  function
compared to the case of a homogeneous medium.

\debproof
We review the results that can be found in 
\cite{andrews,dawson84,furutsu72,furutsu73,garniers3} and put them in a 
convenient form for the derivation.
If the covariance function $\gamma_0$ can be expanded as (\ref{eq:expandgamma0}),
then  the equation for the Fourier transform of the fourth-order moment
can be simplified in the spot-dancing regime $\epsilon  \to 0$ as:
\begin{equation}
 \frac{\partial \hat{\mu}_\br}{\partial z} + \frac{i}{k_o} \big( \bxi_1\cdot \bzeta_1+   \bxi_2\cdot \bzeta_2\big)\hat{\mu}_\br
=
\frac{k_o^2 \bar{\gamma}_2}{2} 
\Delta_{\bxi_1} 
\hat{\mu}_\br.
\end{equation}
This equation can be solved (by a Fourier transform in $\bxi_1$):
\begin{eqnarray}
\nonumber    \hspace*{-2cm}
\hat{\mu}_\br(\bxi_1,\bxi_2,\bzeta_1,\bzeta_2,z) &=&
\int 
\hat{\mu}_\br(\bxi_1',\bxi_2,\bzeta_1,\bzeta_2,0)\\
&&\times
 \exp \Big( -i \frac{z}{k_o} (\bxi_1'\cdot \bzeta_1 +\bxi_2\cdot \bzeta_2) \Big)
\psi(  \bxi_1-\bxi_1',\bzeta_1,z) \d \bxi_1'  ,
\end{eqnarray}
with
\begin{equation}
\psi(  \bxi ,\bzeta_1,z) =\frac{1}{(2 \pi k_o^2  \bar{\gamma}_2 z)^{d/2}}
\exp \Big( - \frac{ \bar{\gamma}_2 z^3}{24} |\bzeta_1|^2 - i \frac{ z}{2k_o} \bxi \cdot \bzeta_1
-\frac{1}{2 k_o^2  \bar{\gamma}_2 z} |\bxi|^2 \Big) .
\end{equation}
This gives  an explicit expression for the fourth-order moment which is what we need to analyze the 
speckle imaging approach considered here. 
As shown in \cite{garniers3}, it is  in fact  possible to compute all the moments 
in the spot-dancing regime and to identify the statistical distribution of the 
transmitted field $E_\br(\bx)$.   We have   in distribution 
\begin{equation}
\hat{E}_\br (\bk) =  \hat{U}_\br \Big( \bk + \frac{k_o\sqrt{  \bar{\gamma}_2}}{2} \bW_{\ell} \Big) \exp 
\Big( - \frac{i}{2k_o} \int_0^{\ell} \big| \bk+ \frac{k_o \sqrt{  \bar{\gamma}_2}}{2} \bW_{z} \big|^2 \d z \Big)  ,
\end{equation}
from which
Proposition \ref{pro:spo} follows. 
\finproofo

\begin{remark}\label{remark3}
To be complete, we can add that it is quite easy to reconstruct the incident field profile $U$ under the natural assumption that
the camera is in the far field (i.e. ${\ell}$ is larger than the Rayleigh length $k_o r_U^2$ where $r_U$ is the radius of the mask).
Indeed, (\ref{eq:fieldfar1}) and (\ref{eq:meanintspot1})  show that the transmitted intensity $| {E}_\br (\bx)|^2$
is equal to $|\hat{U}_{\br+\bX_{\ell}}(k_o \bx / {\ell})|^2$  (up to a multiplicative constant).
From the modulus of the Fourier transform of $U_{\br+\bX_{\ell}}(\bx)$ and from its  phase (assumed to be known, for instance, zero)  
it is possible to reconstruct the incident field profile by a phase-retrieval algorithm \cite{fienup}.
Note,  however, that for a large window  the displacement $\bX_{\ell}$ may vary over the image.   
\end{remark}

\section{Summary and Concluding Remarks}
\label{sec:con}%

We have considered an algorithm for imaging of a moving object based on speckle statistics.
The scheme is as introduced in \cite{webb14} and the basic quantity computed is
the measured or empirical intensity covariance over incident position 
 \ba
\nonumber
C_{\br,\br'} &=& \frac{1}{|A_o|}
\int_{A_o} |E_\br ( \bx_0 )|^2 |E_{\br'}( \bx_0 )|^2 \d \bx_0
\\
& & \hbox{}
- \Big( \frac{1}{|A_o|}
\int_{A_o} |E_\br ( \bx_0 )|^2  \d \bx_0 \Big)
\Big( \frac{1}{|A_o|}
\int_{A_o} |E_{\br'} ( \bx_0 )|^2  \d \bx_0 \Big) ,
\label{def:intcor2}
\ea
where $A_o$ is the spatial support of the  camera and $\br,\br'$ are incident 
positions, see Figure~\ref{fig:1}. 
  The conjecture of \cite{webb14}
is  that  
\begin{equation}
\label{eq:pred2}
C_{\br,\br'}  \approx \Big| 
\int_{\RR^d}
|\hat{U}(\bk)|^2 \exp \big( i\bk \cdot ( \br'-\br) \big)  \d\bk
 \Big|^2  \propto  \left| (U \star \overline{U})( \br - \br') \right|^2  ,
\end{equation}
where $\star$ stands for convolution, 
so that the mask $U$ can be recovered via a phase retrieval step.
The interesting consequence of such a result is that precise information
about the shape of the mask is hidden in the complex speckle pattern, moreover,
that the expression for the empirical intensity covariance does not depend on the  
properties of the complex section and the associated character of the scattering process. 
The argument in \cite{webb14} is based on a strong scattering assumption and an associated 
zero-mean circular Gaussian assumption for the transmitted wave field.  

Here we have presented an analysis of this problem with a view toward
identifying the precise scaling regime where the beautiful relation  \eqref{eq:pred2} 
as set forth in \cite{webb14}  can be mathematically justified when modeling
the complex section as shown in Figure  \ref{fig:1} as a random medium, moreover,
when we consider scalar harmonic wave propagation, as a model for narrow band
optics. 
  
To set the stage for our discussion  let us consider that the random medium 
fluctuations in (\ref{eq:wave0}) have mean zero and covariance of the form
\ban
\EE \big[ \mu (\bx,z) \mu (\bx',z') ] = \sigma^2 {\cal C}_{\mu}\Big( \frac{\bx-\bx'}{\ell_c} , \frac{z-z'}{\ell_c}\Big), 
\ean
with ${\cal C}_{\mu}$ a normalized function (such that ${\cal C}_{\mu}({\bf 0})=1$ and the radius of ${\cal C}_{\mu}$ is of order one).
In this model $\sigma^2$ is the variance of the relative random fluctuations of the medium
and $\ell_c$ is the coherence  length. We also let
\ban
   \gamma_0(\bx-\bx')   =  \int_{-\infty}^{\infty}  \EE \big[ \mu (\bx,z) \mu (\bx',z+z') ]  \, \d z' ,
\ean
which is the lateral spectrum of the driving Brownian motion in  the 
It\^o-Schr\"odinger equation in \eqref{eq:IS}. 

Some central parameters associated with this formulation are then
(i) the central wavelength $\lambda_o=2\pi c_0/k_o$, (ii) the medium coherence length $\ell_c$,
(iii) the relative magnitude of the medium fluctuations $\sigma$, (iii) the radius of the camera  ${r}_A$,
(iv) the size ${r}_U$ of the mask $U$, (v) the distance from the mask to the camera ${\ell}$
corresponding to the thickness of the random section.

The main scaling regime we have considered is the scaling regime leading to the
It\^o-Schr\"odinger equation in \eqref{eq:IS}, or the white-noise paraxial model, corresponding to 
  \ban
  \lambda_o=2\pi/k_o \ll \ell_c  \ll {\ell}     .
\ean
Then, we have considered two subregimes of propagation which essentially
are the two canonical scaling regimes in the  white-noise paraxial model:
(a) the {\it scintillation regime}  corresponding to ${r}_U \gg \ell_c$,
(b) the {\it spot-dancing regime} corresponding  to ${r}_U \ll \ell_c$.

In the spot-dancing regime,  the wave  
intensity pattern is as in the homogeneous case, however, modified
by a random lateral shift in the profile. In fact, in this case the formula 
\eqref{eq:pred2} is not valid, however, the mask can still be recovered,
albeit with a different approach corresponding to the one one would 
have used in a homogeneous medium.  

In the scintillation regime, 
the transmitted  wave forms a speckle pattern with rapid fluctuations
of the intensity. 
In order to discuss the scintillation regime let us introduce two parameters.
First,  the characteristic size of the speckle fluctuations or speckle radius at range ${\ell}$ is  
\ban
\rho_{\ell} =   \frac{\ell_c^{1/2}}{ \sigma k_o {\ell}^{1/2} } .
\ean
 The other fundamental parameter associated with the scintillation regime 
is the beam spreading  width at range ${\ell}$ which is
\ban
   {\mathcal A}_{\ell} =  \frac{\sigma   {\ell}^{3/2}}{\ell_c^{1/2}}  =   \frac{k_o {\ell}}{\rho_{\ell}} . 
\ean 
In order to have a high signal-to-noise ratio so that the empirical intensity covariance
function is close to its expectation we assume
\ban
    \rho_{\ell} \ll {r}_A .
\ean
We remark that if the camera is associated with finite-sized elements,
of size $\rho_o$, then we assume that $ \rho_o =O (\rho_{\ell})$ to retain 
a high signal-to-noise ratio  (the effects of having  finite-sized elements is analyzed 
in detail above). \\
We then  arrive at the asymptotic description  in \eqref{eq:c} for the empirical
intensity covariance. This expression involves the
 medium (second-order) statistics and the mask function $U$. 
It can form the basis for an estimation procedure for the mask 
and we remark that  it holds true whatever the magnitude of the 
scattering  mean free path $\ell_{\rm sca}$ is relative to the range ${\ell}$, with
$\ell_{\rm sca}$ given in \eqref{def:lsca:parax} and which  corresponds to 
\ban
 \ell_{\rm sca}   = {\ell} \Big( \frac{\rho_{\ell}}{\ell_c} \Big)^2 ,
\ean
so that the regime of long-range propagation corresponds to $\rho_{\ell} \ll \ell_c$. 
Upon some last scaling  assumptions we arrive exactly
at the description in \eqref{eq:pred2}. Specifically assume
(i) relatively large spreading so that $|\br|, {r}_A \ll {\mathcal A}_{\ell}$
(ii)  long-range propagation so that $\ell_{\rm sca}  \ll {\ell}$ and (iii) smooth medium fluctuations
so that \eqref{eq:expandgamma0} is valid.  These are the last stepping stones toward the 
formula \eqref{eq:pred2}.    
 
Let us next comment on an informal interpretation of the above  result.
Let $G_{\ell}(\bx,\br)$ be  the Green's function  over the section
$z\in (0,{\ell})$ for a source point at $(\br,0)$ and  an observation point at  $(\bx,{\ell})$.
Then we have for the transmitted field:
\ban
   E_{\br}(\bx)   =      \int_{\RR^d} U(\by-\br) G_{\ell}(\bx, \by) \d \by   .
 \ean
Let us first consider  the   
 field covariance function with respect to shift vector $\br$:
 \ban
 {\cal D}_{{\br},{\br'}}(\bx_0,\bx_0)      
 &=& \EE \big[ E_\br ( \bx_0 ) \overline{E_{\br'}( \bx_0 )}\big]
 \ean
  which,  
 making use of reciprocity,
 can be  expressed as 
\ban
 {\cal D}_{{\br},{\br'}}(\bx_0,\bx_0)       &=&  \iint_{\RR^{2d}} 
   U(\by-\br) \overline{U(\tilde{\by}-\br')} 
   \EE\left[ G_{\ell}(\bx_0, \by)   \overline{G_{\ell}(\bx_0, \tilde{\by}) }  \right] 
    \d \by  \d \tilde{\by}  .
\ean
In the strongly scattering regime and under the assumption that the speckle radius $\rho_{\ell}$ is much smaller than ${r}_U$,
the covariance $ \EE [ G_{\ell}(\bx_0, \by)   \overline{G_{\ell}(\bx_0, \tilde{\by}) } ] $ is approximately delta-correlated in $\by-\tilde{\by}$
and is proportional to an envelope with beam width ${\cal A}_{\ell}$, so that
we get
\ban
 &&   {\cal D}_{{\br},{\br'}}(\bx_0,\bx_0)         \propto   \int_{\RR^{d}} 
   U(\by- \br) \overline{U( {\by} -\br')}  
    s\Big(   \frac{\by-\bx_0 }{{\mathcal A}_{\ell}} \Big) 
    \d \by      ,
 \ean
for $s$ a normalized envelope function with unit width and unit amplitude.
Under the assumption that ${r}_U,{r}_A$ and the camera center point have small 
magnitude relative to ${\mathcal A}_{\ell}$ we get   
\ban
 &&   {\cal D}_{{\br},{\br'}}(\bx_0,\bx_0)        \propto    (U \star \overline{U})( \br' - \br)  .
  \ean 
 We next  have for the speckle covariance function with respect to the shift vector   
 \ban
 && \hspace*{-2.7cm} {\cal C}_{{\br},{\br'}}(\bx_0,\bx_0)   =\EE\big[ |E_\br ( \bx_0 )|^2 |E_{\br'}( \bx_0 )|^2 \big]
-
\EE\big[ |E_\br ( \bx_0 )|^2 \big]\EE\big[|E_{\br'}( \bx_0 )|^2 \big]
 \\&=&   \EE\bigg[
          \int_{\RR^d} U(\by-\br) G_{\ell}(\bx_0, \by) \d \by   \int_{\RR^d} U(\by-\br') G_{\ell}(\bx_0, \by) \d \by
         \\ &  & \hbox{}  \times     
          \int_{\RR^d}   \overline{U(\by-\br) G_{\ell}(\bx_0, \by) }\d \by  
             \int_{\RR^d}   \overline{U(\by-\br') G_{\ell}(\bx_0, \by) }\d \by  \bigg]   \\
             &  &
              - {\cal D}_{{\br},{\br}}(\bx_0,\bx_0)  \overline{ {\cal D}_{{\br'},{\br'}}(\bx_0,\bx_0)  }
               \\   &=& 
             \iint_{\RR^{4d}}   \left( U(\by_1-\br)  \overline{U(\by_2-\br')} \right)
              \left( \overline{U(\by_3-\br)}  {U(\by_4-\br')} \right)   
              \\ &  & \hbox{}  \times
              \EE\left[  G_{\ell}(\bx_0, \by_1) \overline{G_{\ell}(\bx_0, \by_2)}   
               \overline{G_{\ell}(\bx_0, \by_3)}  {G_{\ell}(\bx_0, \by_4)}  \right] 
                  \d \by_1  \d \by_2  \d \by_3  \d \by_4 \\
                  & & \hbox{}
                  - {\cal D}_{{\br},{\br}}(\bx_0,\bx_0)  \overline{ {\cal D}_{{\br'},{\br'}}(\bx_0,\bx_0)  }
                    \\  &=&
                   {\cal D}_{{\br},{\br'}}(\bx_0,\bx_0)        \overline{ {\cal D}_{{\br},{\br'}}(\bx_0,\bx_0)} ,
 \ean
 where we  have used a Gaussian summation rule (Isserlis formula) which states that for 
  four   jointly  complex circularly symmetric Gaussian random variables, 
  $Z_j, j=1,\ldots,4$,        we have 
\begin{eqnarray}
\label{eq:gaussrule}
 \EE \big[  Z_1 \overline{Z_2  Z_3} Z_4\big]  =
 \EE\big[  Z_1  \overline{Z_2} \big] \EE \big[\overline{ Z_3 } Z_4\big] 
+
 \EE\big[  {Z_1  } \overline{Z_3} \big] \EE \big[ \overline{ Z_2 } Z_4\big] .
\end{eqnarray}
  We then arrive at     
   \ban
   &&   {\cal C}_{{\br},{\br'}}(\bx_0,\bx_0)   \propto  \left|  ( U \star \overline{U})( \br - \br')  \right|^2 ,
    \ean
 which is \eqref{eq:pred2}. 
 We   comment here that it is clear from the above argument that in this version
of speckle imaging the so-called memory effect for the speckle pattern, which
is important in other modalities  of speckle imaging  \cite{freund88,vellekoop07,vellekoop08},  is not important. What 
is important here is a small speckle radius and a large spreading of the field. 
Moreover, in this formal argument we made use of a Gaussian assumption which made it possible
to factor a fourth moment in terms of second moments. That this is valid in the considered
regime is a deep result of waves in random media which was recently developed in 
\cite{garniers4}.  
  Note also that the above argument  shows how a similar mask imaging procedure
can be constructed when  we have access to the wave field itself:
it is then possible to estimate the field covariance function with respect to shift vector 
and the Gaussian property is not needed. 
   
Finally,  in Remarks \ref{remark2} and \ref{remark3} we discuss
how, under various circumstances about the random medium,  
the image may be subject to  blurring and geometric distortion operators.
In practice some amount of  both of these effects will be present. For instance, in
the context of turbulence mitigation for  propagation
through the   atmosphere, they need to be corrected for.
We refer to \cite{gilles0,gilles1,bert}
for  frameworks that aim  at mitigating such effects  
 where in particular a physical model, the so-called ``fried kernel'',   is 
 partly and successfully being used. 
Here, we have developed the theory for how such distortion operators
can be modeled in the context of speckle imaging.
{Indeed in this paper we were able to address 
separately the two canonical scaling regimes in the  white-noise paraxial model:
the  scintillation regime   corresponding to ${r}_U \gg \ell_c$ and 
the spot-dancing regime  corresponding  to ${r}_U \ll \ell_c$.
The intermediate regime, when ${r}_U \sim \ell_c$, cannot be addressed via the asymptotic techniques used in our paper.
We may expect  that it should produce  a mixture 
of the two canonical scaling regimes, which would result
in a more  challenging situation from the inverse problems point of view. In particular, we anticipate that the 
intensity covariance function should then not be statistically stable.}

\section*{Acknowledgments}  
This research is supported in part by 
 AFOSR grant  FA9550-18-1-0217, NSF  grant 1616954, 
 Centre Cournot, Fondation Cournot, and 
 Universit\'e Paris Saclay (chaire D'Alembert).

\bigskip


\medskip
Received xxxx 20xx; revised xxxx 20xx.
\medskip

\end{document}